\documentclass[leqno,12pt]{article}

\usepackage{amssymb}
\usepackage{amsthm}
\usepackage{euscript}
\usepackage{fancyhdr}
\usepackage{graphicx}

\evensidemargin\oddsidemargin

\begin{document}

\baselineskip=18pt 

\renewcommand{\theequation}{\thesection.\arabic{equation}}
\newtheorem{theorem}{Theorem}[section]
\newtheorem{lemma}[theorem]{Lemma}
\newtheorem{proposition}[theorem]{Proposition}
\newtheorem{corollary}[theorem]{Corollary}
\newtheorem{remark}[theorem]{Remark}
\newtheorem{fact}[theorem]{Fact}
\newtheorem{problem}[theorem]{Problem}
\newtheorem{corollaire}[theorem]{Corollaire}
\newcommand{\eqnsection}{
\renewcommand{\theequation}{\thesection.\arabic{equation}}
    \makeatletter
    \csname  @addtoreset\endcsname{equation}{section}
    \makeatother}
\eqnsection
\newcommand{\cE}{\mathcal{E}}
\newcommand{\cF}{\mathcal{F}}
\newcommand{\bR}{\mathbb{R}}
\newcommand{\cS}{\mathcal{S}}
\newcommand{\cT}{\mathcal{T}}
\newcommand{\ho}{{\rm H\"o}}
\newcommand{\unif}{{\rm \bf {UNIF}}}
\newcommand{\Prob}{{\rm P}}
\newcommand{\GW}{{\rm GW}}

\vglue15pt
\centerline{\LARGE The uniform measure on a Galton--Watson tree  }

\bigskip
\centerline{\LARGE without the XlogX condition}

\bigskip\bigskip
\centerline{Elie A\"id\'ekon\footnote{ Department of Mathematics and Computer science, Eindhoven University of Technology, P.O. Box 513, 5600 MB Eindhoven, The Netherlands. email: elie.aidekon@gmail.com}}
\bigskip

\bigskip
\bigskip

{\leftskip=2truecm \rightskip=2truecm \baselineskip=15pt \small

\noindent{\slshape\bfseries Summary.} We consider a Galton--Watson tree with offspring distribution $\nu$ of finite mean. The uniform measure on the boundary of the tree is obtained by putting mass $1$ on each vertex of the $n$-th generation and taking the limit $n\to \infty$. In the case $E[\nu\ln(\nu)]<\infty$, this measure has been well studied, and it is known that the Hausdorff dimension of the measure is equal to $\ln(m)$ (\cite{hawkes}, \cite{lpp95}). When $E[\nu \ln(\nu)]=\infty$, we show that the dimension drops to $0$. This answers a question of Lyons, Pemantle and Peres \cite{LyPemPer97} }.

\bigskip
\bigskip

{\leftskip=2truecm \rightskip=2truecm \baselineskip=15pt \small

\noindent{\slshape\bfseries R\'esum\'e.} Nous consid\'erons un arbre de Galton--Watson dont le nombre d'enfants $\nu$ a une moyenne finie. La mesure uniforme sur la fronti\`ere de l'arbre s'obtient en chargeant chaque sommet de la $n$-i\`eme g\'en\'eration avec une masse $1$, puis en prenant la limite $n\to\infty$. Dans le cas $E[\nu\ln(\nu)]<\infty$, cette mesure a \'et\'e tr\`es \'etudi\'ee, et l'on sait que la dimension de Hausdorff de la mesure est \'egale \`a $\ln(m)$ (\cite{hawkes}, \cite{lpp95}). Lorsque $E[\nu \ln(\nu)]=\infty$, nous montrons que la dimension est $0$. Cela r\'epond \`a une question pos\'ee par Lyons, Pemantle et Peres \cite{LyPemPer97} }.



\section{Introduction}

Let $\cT$ be a Galton--Watson tree of root $e$, associated to the offspring distribution $q:=(q_k,k\ge 0)$. We denote by $\GW$ the distribution of $\cT$ on the space of rooted trees, and $\nu$ a generic random variable on $\mathbb{N}$ with distribution $q$. We suppose that $q_0=0$ and $m:=\sum_{k\ge 0} k q_k \in (1,\infty)$: the tree has no leaf (hence survives forever) and is not degenerate. For any vertex $u$, we write $|u|$ for the height of vertex $u$ ($|e|=0$), $\nu(u)$ for the number of children of $u$, and $Z_n$ is the population at height $n$. We define $S(\cT)$ as the set of all infinite self-avoiding paths of $\cT$ starting from the root and we define a metric on $S(\cT)$ by $d(r,r'):=e^{-|r \land r'|}$ where $r\land r'$ is the highest vertex belonging to $r$ and $r'$. The space $S(\cT)$ is called boundary of the tree, and elements of $S(\cT)$ are called rays. \\

When $E[\nu\ln(\nu)]<\infty$, it is well-known that the martingale $m^{-n}Z_n$ converges in $L^1$ and almost surely to a positive limit (\cite{kestenstigum}). Seneta \cite{seneta} and Heyde \cite{heyde} proved that in the general case (i.e allowing $E[\nu\ln(\nu)]$ to be infinite), there exist constants $(c_n)_{n\ge 0}$ such that
\begin{itemize}
\item[(a)] $W_{\infty}:=\lim_{n\to \infty}{Z_n \over c_n}$ exists a.s.
\item[(b)] $W_{\infty}>0$ a.s.
\item[(c)] $\lim_{n\to \infty} {c_{n+1} \over c_n} = {m}$.
\end{itemize}

\noindent In particular, for each vertex $u\in\cT$, if $Z_k(u)$ stands for the number of descendants $v$ of $u$ such that $|v|=|u|+k$, we can define
    $$
    W_{\infty}(u):=\lim_{k\to\infty} {Z_k(u)\over c_k}
    $$

\noindent and we notice that $m^{-n}\sum_{|u|=n} W_{\infty}(u)=W_{\infty}(e)$.

\bigskip

\noindent {\bf Definition}. The uniform measure (also called branching measure) is the unique Borel measure on $S(\cT)$ such that
    $$
    \unif(\{r\in S(\cT), \, r_n=u\}) := {m^{-n} W_{\infty}(u) \over W_{\infty}(e)}
    $$
for any integer $n$ and any vertex $u$ of height $n$. \\

\noindent We observe that, for any vertex $u$ of height $n$,
    $$
    \unif(\{r\in S(\cT), \, r_n=u\}) = \lim_{k\to \infty} {Z_{k}(u) \over Z_{n+k}}\,.
    $$
Therefore the uniform measure can be seen informally as the probability distribution of a ray taken uniformly in the boundary. This paper is interested in the Hausdorff dimension of $\unif$, defined by
    $$
    \dim(\unif) := \min\{\dim(E),\, \unif(E)= 1\}
    $$
where the minimum is taken over all subsets $E\subset S(\cT)$ and $\dim(E)$ is the Hausdorff dimension of set $E$. The case $E[\nu\ln(\nu)]<\infty$ has been well studied. In \cite{hawkes} and \cite{lpp95}, it is shown that $\dim(\unif) = \ln(m)$ almost surely. A description of the multifractal spectrum is available in \cite{Liu01},\cite{lpp95},\cite{MortersShieh},\cite{shiehtaylor}. The case $E[\nu \ln(\nu)]=\infty$ presented as Question 3.1 in \cite{LyPemPer97} was left open. This case is proved to display an extreme behaviour.

\begin{theorem}
\label{hausdorff}
If $E[\nu \ln(\nu)]=\infty$, then $\dim(\unif)=0$ for $\GW$-a.e tree $\cT$.
\end{theorem}

\bigskip

The drop in the dimension comes from bursts of offspring at some places of the tree $\cT$. Namely, for $\unif$-a.e. ray $r$, the number of children of $r_n$ will be greater than $(m-o(1))^n$ for infinitely many $n$. To prove it, we work with a particular measure $Q$, under which the distribution of the numbers of children of a uniformly chosen ray is more tractable. Section \ref{sect_spine} contains the description of the new measure in terms of a spine decomposition. Then we prove Theorem \ref{hausdorff} in Section \ref{sect_proof}.

\section{A spine decomposition}\label{sect_spine}

For $k\ge 1$ and $s\in (0,1)$, we call $\phi_k(s)$ the probability generating function of $Z_k$
    $$
    \phi_k(s) := E[s^{Z_k}]\,.
    $$

\noindent We denote by $\phi^{-1}_k(s)$ the inverse map on $(0,1)$ and we let $s\in (0,1)$. Then  $M_n:=\phi_n^{-1}(s)^{Z_n}$ defines a martingale and converges in $L^1$ to some $M_{\infty}>0$ a.s (\cite{heyde}). Therefore we can take in (a)
    $$
    c_n := {-1 \over  \ln(\phi_n^{-1}(s))}
    $$
\noindent which we will do from now on. Hence we can rewrite equivalently $M_n= e^{-Z_n/c_n}$ and $M_{\infty}=e^{-W_{\infty}(e)}$. For any vertex $u$ at generation $n$, we define similarly
    $$
    M_{\infty}(u):=  {1\over \phi_n^{-1}(s)}e^{-m^{-n}W_{\infty}(u)} = e^{1/c_n}e^{-m^{-n}W_{\infty}(u)}
    $$
which is the limit of the martingale $M_{k}(u):=e^{1/c_n}e^{ - Z_k(u)/c_{n+k}}$. In \cite{lynch}, Lynch introduces the so-called derivative martingale
    $$
    \partial M_n:= e^{1/c_n}{Z_n\over \phi'_n(\phi_n^{-1}(s))}M_n
    $$

\noindent and shows that the derivative martingale also converges almost surely and in $L^1$ ($\partial M_n$ is in fact bounded). Moreover the limit $\partial M_{\infty}$ is positive almost surely. We deduce that the ratio $\phi'_n(\phi_n^{-1}(s))/c_n$ converges to some positive constant. In particular, it follows from (c) that
    \begin{equation}\label{phi'}
    \lim_{n\to \infty} {\phi'_{n+1}(\phi_{n+1}^{-1}(s)) \over \phi'_n(\phi_n^{-1}(s))} ={m}.
    \end{equation}

\bigskip
We are interested in the probability measure $Q$ on the space of rooted trees defined by
    $$
    {dQ\over d\GW}:= \partial M_{\infty}\,.
    $$

\noindent Let us describe this change of measure. We call a marked tree a couple $(T,\, r)$ where $T$ is a rooted tree and $r$ a ray of the tree $T$. Let $(\mathbb{T},\,\xi)$ be a random variable in the space of all marked trees (equipped with some probability $\Prob(\cdot)$), whose distribution is given by the following rules. Conditionally on the tree up to level $k$ and on the location of the ray at level $k$, (which we denote respectively by $\mathbb{T}_k$ and $\xi_k$),
\begin{itemize}
    \item the number of children of the vertices at generation $k$ are independent
	\item the vertex $\xi_k$ has a number $\nu(\xi_k)$ of children such that for any $\ell$
	\begin{equation}\label{offspring_spine}
	\Prob(\nu(\xi_k)=\ell)=  \tilde q_{\ell}^s := q_{\ell}\, \ell \exp\left(-{\ell-1 \over c_{k+1}}\right) {\phi_{k}'(\phi_{k}^{-1}(s))\over \phi_{k+1}'(\phi_{k+1}^{-1}(s)) }
	\end{equation}
	\item the number of children of a vertex $u\neq \xi_k$ at generation $k$  verifies for any $\ell$
	\begin{equation}\label{offspring_notspine}
	\Prob(\nu(u)=\ell)=\tilde q_{\ell}:= q_{\ell} \,e^{1/c_k}\exp\left(-{\ell \over c_{k+1}}\right)
	\end{equation}
	\item the vertex $\xi_{k+1}$ is chosen uniformly among the children of $\xi_k$
\end{itemize}

\noindent As often in the literature, we will call the ray $\xi$ the spine. We refer to \cite{LyoPemPer}, \cite{lyons97} for motivation on spine decompositions. In our case, we can see $\mathbb{T}$ as a Galton-Watson tree in varying environment and with immigration. The fact that (\ref{offspring_spine}) and (\ref{offspring_notspine}) define probabilities come from the equations (remember that by definition $e^{-1/c_k}=\phi_k^{-1}(s)$)
    \begin{eqnarray*}
    E\left[ (\phi_{k+1}^{-1}(s))^{\nu} \right]& =& \phi_k^{-1}(s)\,, \\
    E\left[\nu (\phi_{k+1}^{-1}(s))^{\nu-1}\right] & = & {\phi_{k+1}'(\phi_{k+1}^{-1}(s))\over \phi_{k}'(\phi_{k}^{-1}(s)) }\,.
    \end{eqnarray*}

\bigskip

\noindent We mention that in \cite{LyoPemPer}, a similar decomposition was presented using the martingale ${Z_n \over m^n}$. In this case, the offspring distribution of the spine is the size-biased distribution $({\ell q_{\ell}\over m})_{\ell\ge 0}$ whereas the other particles generate offspring according to the original distribution $q$. In particular, the offspring distributions do not depend on the generation. When $E[\nu\ln(\nu)]<\infty$, the process, which is a Galton--Watson process with immigration, has a distribution equivalent to $\GW$. It is no longer true  when $E[\nu \ln(\nu)]=\infty$, in which case the spine can give birth to a super-exponential number of children.

\bigskip

\begin{proposition} \label{spine}
Under $Q$, the tree $\cT$ has the distribution of $\mathbb{T}$. Besides, for $\Prob$-almost every tree $\mathbb{T}$, the distribution of $\xi$ conditionally on $\mathbb{T}$ is the uniform measure $\unif$.
\end{proposition}

\noindent {\it Proof}. For any tree $T$, we define ${T}_n$ the tree ${T}$ obtained by keeping only the $n$-first generations. Let $T$ be a tree. We will prove by induction that, for any integer $n$  and any vertex $u$ at generation $n$,
    \begin{equation}\label{induction}
    \Prob(\mathbb{T}_{ n} =T_{ n},\,\xi_{n}=u) = {\partial M_n \over Z_n} \GW(\cT_n =T_{ n})\,.
    \end{equation}

\noindent For $n=0$, it is straightforward since $\mathbb{T}_0$ and $\cT_0$ are reduced to the root. We suppose that this is true for $n-1$, and we prove it for $n$. Let ${\buildrel \leftarrow \over u}$ denote the parent of $u$, and, for any vertex $v$ at height $n-1$, let $k(v)$ denote the number of children of $v$ in the tree $T$. We have
    \begin{eqnarray*}
    && \Prob(\mathbb{T}_{ n} =    T_{ n},\,\xi_{n}=u\,|\, \mathbb{T}_{n-1}=T_{n-1},\, \xi_{n-1} = {\buildrel \leftarrow \over u})\\
    &=&
    {1\over k({\buildrel \leftarrow \over u})}{\tilde q^s_{k({\buildrel \leftarrow \over u})} \over \tilde q_{k({\buildrel \leftarrow \over u})}} \prod_{|v|=n-1}\tilde q_{k(v)}\\
    &=&
    {e^{1/c_n} \over e^{1/c_{n-1}}}{\phi_{n-1}'(\phi_{n-1}^{-1}(s))\over \phi_{n}'(\phi_{n}^{-1}(s)) }{ e^{Z_{n-1}/c_{n-1}} \over e^{Z_n / c_n}} \prod_{|v|=n-1}q_{k(v)}\\
    &=&
    {e^{1/c_n} \over e^{1/c_{n-1}}}{\phi_{n-1}'(\phi_{n-1}^{-1}(s))\over \phi_{n}'(\phi_{n}^{-1}(s)) }{ e^{Z_{n-1}/c_{n-1}} \over e^{Z_n / c_n}}
    \GW(\cT_n = T_n \,|\, \cT_{n-1}=T_{n-1}) \,.
    \end{eqnarray*}

\noindent We use the induction assumption to get
    \begin{eqnarray*}
    \Prob(\mathbb{T}_{ n} =    T_{ n},\,\xi_{n}=u) = e^{1/c_n}{1 \over \phi_{n}'(\phi_{n}^{-1}(s)) }e^{-Z_n \over c_n}
    \GW(\cT_n = T_n)
    \end{eqnarray*}

\noindent which proves (\ref{induction}). Summing over the $n$-th generation of $T$ gives
$$
\Prob(\mathbb{T}_{ n} =T_{ n}) = \partial M_n \GW(\cT_n=T_n) = Q(\cT =T_{ n})\,.
$$

\noindent This computation also shows that $\Prob(\xi_n= u\,|\, \mathbb{T}_{ n})= 1/Z_n$ which implies that $\xi$ is uniformly distributed on the boundary $S(\mathbb{T})$. $\Box$

\bigskip

\noindent {\bf Remark A}. For $u$ a vertex of $\mathbb{T}$ at generation $n$, call $\mathbb{T}(u)$ the subtree rooted at $u$. A similar computation shows that if $u\notin \xi$, then the distribution $\Prob_u$ of $\mathbb{T}(u)$ (conditionally on $\mathbb{T}_n$ and on $\xi_{n}$) verifies
    $$
    {d\Prob_u \over d\GW} = M_{\infty}(u)\,.
    $$

\section{Proof of Theorem \ref{hausdorff}} \label{sect_proof}

The following proposition shows that in the tree $\cT$, there exist infinitely many times when the ball $\{r \in S(\cT)\,:\, r_n=\xi_n\}$ has a 'big' weight.
\begin{proposition}\label{holdprob}
Suppose that $E[\nu\ln(\nu)]=\infty$. Then we have $\Prob$-a.s.
    $$
    \limsup_{n\to\infty} {1\over n} \ln(W_{\infty}(\xi_n)) = \ln(m)\,.
    $$
\end{proposition}

\noindent {\it Proof}. Let $1<a<b<m$ and $n\ge 0$. We get from (\ref{offspring_spine})
    \begin{eqnarray*}
    \Prob\left(\nu(\xi_n)\in  (a^n,b^n)\right)
    &=&
    {\phi_{n}'(\phi_{n}^{-1}(s))\over \phi_{n+1}'(\phi_{n+1}^{-1}(s)) }E\left[\nu e^{-(\nu-1)/c_{n+1}},\, \nu \in (a^n,b^n) \right]\\
    &\ge&
    {\phi_{n}'(\phi_{n}^{-1}(s))\over \phi_{n+1}'(\phi_{n+1}^{-1}(s)) }e^{-b^n/c_n}E\left[\nu,\, \nu \in (a^n,b^n) \right]
    \,.
\end{eqnarray*}

\noindent From (c) and (\ref{phi'}), we deduce that for $n$ large enough, we have
\begin{eqnarray*}
\Prob\left(\nu(\xi_n) \in (a^n,b^n)\right) \ge   {1\over 2m} E\left[\nu,\, \nu \in (a^n,b^n) \right]\,.
\end{eqnarray*}

\noindent Therefore, under the condition $E[\nu\ln(\nu)]=\infty$, we have
    \begin{equation}\label{inf_sum}
    \sum_{n\ge 0} \Prob(\nu(\xi_n) \in (a^n,b^n)) = \infty\,.
    \end{equation}

\noindent Let $H(\xi_n):=\{u\in \mathbb{T}: u \mbox{ child of } \xi_n,\, u\neq \xi_{n+1}\}$. By Remark A, we have
    \begin{eqnarray*}
     \Prob\left( \sum_{u \in H(\xi_n)} W_{\infty}(u) \le a^n  \,\bigg|\, \mathbb{T}_{n+1},\xi_{n+1}\right)
    =
    E_{\GW}\left[ \prod_{u\in \mathcal H} M_{\infty}(u), \sum_{u \in \mathcal H} W_{\infty}(u) \le a^n \right]_{\mathcal H = H(\xi_n)}\,.
    \end{eqnarray*}

\noindent Since $M_{\infty}(u)\le e^{1/c_n}$ for any $|u|=n$, we get
    \begin{eqnarray*}
    \Prob\left( \sum_{u \in H(\xi_n)} W_{\infty}(u) \le a^n  \,\bigg|\, \mathbb{T}_{n+1},\xi_{n+1}\right)
    \le
    e^{(\nu(\xi_n)-1)/c_n} \GW\left(\sum_{u \in \mathcal{H}} W_{\infty}(u) \le a^n  \right)_{\mathcal{H}=H(\xi_n)}\,.
    \end{eqnarray*}

\noindent Let $(W_{\infty}^i,i\ge 1)$ be independent random variables distributed as $W_{\infty}(e)$ under $\GW$. It follows that on the event $\{\nu(\xi_n) \in (a^n,b^n)\}$, we have
    \begin{eqnarray*}
    \Prob\left( \sum_{u \in H(\xi_n)} W_{\infty}(u) \le a^n  \,\bigg|\, \mathbb{T}_{n+1},\xi_{n+1}\right)
    \le
    e^{(b^n-1)/c_n}\GW\left(\sum_{i=1}^{a^n} W_{\infty}^i \le a^n\right)=:d_n\,.
    \end{eqnarray*}

\noindent We obtain that
    \begin{eqnarray*}
    \Prob\left(\sum_{u \in H(\xi_n)} W_{\infty}(u) > a^n\right)
    &\ge&
    \Prob\left(\sum_{u \in H(\xi_n)} W_{\infty}(u)> a^n,\,\nu(\xi_n) \in (a^n,b^n)\right)\\
    &\ge&
    \Prob\left(\nu(\xi_n) \in (a^n,b^n)\right)(1-d_n)\,.
    \end{eqnarray*}

\noindent By (c), $e^{(b^n-1)/c_n}$ goes to $1$. Furthermore, we know from \cite{seneta} that $E_{\GW}[W_{\infty}(e)]=\infty$, which ensures by the law of large numbers that $d_n$ goes to $0$. By equation (\ref{inf_sum}), we deduce that
    $$
    \sum_{n\ge 0} \Prob\left( \sum_{u \in H(\xi_n)} W_{\infty}(u) >  a^n\right) = \infty\,.
    $$

\noindent We use the Borel-Cantelli lemma to see that $\sum_{u \in H(\xi_n)} W_{\infty}(u) >  a^n$ infinitely often. Since $W_{\infty}(\xi_n) \ge {1\over m} \sum_{u\in H(\xi_n)} W_{\infty}(u)$, we get that $W_\infty(\xi_n) \ge a^n/m$ for infinitely many $n$, $\Prob$-a.s. Hence
    $$
    \limsup_{n\to \infty}{1\over n}\ln(W_\infty(\xi_n)) \ge \ln(a)\,.
    $$

\noindent Let $a$ go to $m$ to have the lower bound. Since $W_\infty(\xi_n) \le \sum_{|u|=n} W_{\infty}(u) = m^n W_{\infty}(e)$, we have $\limsup_{n\to\infty} {1\over n}\ln(W_\infty(\xi_n)) \le \ln(m)$ hence Proposition \ref{holdprob}. $\Box$

\bigskip

We turn to the proof of the theorem.

\bigskip

\noindent{\it Proof of Theorem \ref{hausdorff}}. By Proposition \ref{holdprob}, we have
    $$
    \Prob\left(\limsup_{n\to\infty} {1\over n} \ln(W_{\infty}(\xi_n)) = \ln(m)\right)=1\,.
    $$

\noindent In particular, for $\Prob$-a.e. $\mathbb{T}$,
    $$
    \Prob\left(\limsup_{n\to\infty} {1\over n} \ln(W_{\infty}(\xi_n)) = \ln(m)\,\bigg|\, \mathbb{T}\right) =1\,.
    $$

\noindent By Proposition \ref{spine}, the distribution of $\xi$ given $\mathbb{T}$ is $\unif$. Therefore, for $\Prob$-a.e. $\mathbb{T}$,
    \begin{equation}\label{unif}
    \unif\left(r\in S(\mathbb{T})\,:\,\limsup_{n\to\infty} {1\over n} \ln(W_{\infty}(r_n)) = \ln(m) \right) = 1\,.
    \end{equation}

\noindent Again by Proposition \ref{spine}, the distribution of $\mathbb{T}$ is the one of $\cT$ under $Q$. We deduce that (\ref{unif}) holds for $Q$-a.e. tree $\cT$. Since $Q$ and $\GW$ are equivalent, equation (\ref{unif}) holds for $\GW$-a.e. tree $\cT$. We call the H\"older exponent of $\unif$ at ray $r^*$ the quantity
    $$
    \mbox{H\"o}(\unif)(r^*) := \liminf_{n\to \infty}{-1\over n} \ln\left(\unif(\{r\in S(\cT)\,:\, r_n =r^*_n\})\right)\,.
    $$

\noindent By definition of $\unif$, we can rewrite it
    $$
    \mbox{H\"o}(\unif)(r^*) = \liminf_{n\to \infty}{-1\over n} \ln\left( m^{-n} W_{\infty}(r_n^*)/W_{\infty}(e)\right)\,.
    $$

\noindent Therefore, for $\unif$-a.e. ray $r$, $\mbox{H\"o}(\unif)(r)=0$. By Theorem 14.15 of \cite{LyPe} (or $\mathsection$ 14 of \cite{billingsley_ergo}), it implies that $\dim(\unif)=0$ $\GW$-almost surely. $\Box$

\bigskip

\noindent {\it Acknowledgements}. The author thanks Russell Lyons for useful comments on the work. This work was supported in part by the Netherlands Organisation for Scientific Research (NWO).

\end{document}